\input amstex
\documentstyle{amams} 
\document
\annalsline{155}{2002}
\received{September 15, 2000}
\startingpage{209}
\font\tenrm=cmr10
\def\therosteritem#1{(#1)}
\catcode`\@=11
\font\twelvemsb=msbm10 scaled 1100

\font\ninemsb=msbm10 scaled 800
\newfam\msbfam
\textfont\msbfam=\twelvemsb  \scriptfont\msbfam=\ninemsb
  \scriptscriptfont\msbfam=\ninemsb
\def\msb@{\hexnumber@\msbfam}
\def\Bbb{\relax\ifmmode\let\next\Bbb@\else
 \def\next{\errmessage{Use \string\Bbb\space only in math
mode}}\fi\next}
\def\Bbb@#1{{\Bbb@@{#1}}}
\def\Bbb@@#1{\fam\msbfam#1}
\catcode`\@=12

 \catcode`\@=11
\font\twelveeuf=eufm10 scaled 1100
\font\teneuf=eufm10
\font\nineeuf=eufm7 scaled 1100
\newfam\euffam
\textfont\euffam=\twelveeuf  \scriptfont\euffam=\teneuf
  \scriptscriptfont\euffam=\nineeuf
\def\euf@{\hexnumber@\euffam}
\def\frak{\relax\ifmmode\let\next\frak@\else
 \def\next{\errmessage{Use \string\frak\space only in math
mode}}\fi\next}
\def\frak@#1{{\frak@@{#1}}}
\def\frak@@#1{\fam\euffam#1}
\catcode`\@=12
\font\titr=cmr10 scaled \magstep3
 \font\esmi= cmmi10  scaled \magstep3
\title{Embeddedness of minimal surfaces with\\ total boundary curvature  at most  \hbox{\titr 4}\hbox{\esmi
\char26}}  
\shorttitle{Embeddedness of minimal surfaces} 
\acknowledgements{The authors were supported in part by:
       (Ekholm) Stiftelsen f\"or internationalisering av forskning
       och h\"ogre utbildning,
       (White) National Science Foundation Grant DMS-9803493 and 
       the Guggenheim Foundation, 
       and (Wienholtz) the Deutsche Forschungsgemeinschaft.\hfill\break
\hglue26pt 1991 {\it Mathematics Subject
Classification}. Primary 53A10, secondary 49F10.\hfill\break 
\hglue26pt {\it Key words and phrases}. Minimal surfaces, total
curvature.}
 \twoauthors{Tobias Ekholm, Brian White,}{Daniel Wienholtz}
   \institutions{Uppsala University, 752 06 Uppsala, Sweden\\
 {\eightpoint {\it E-mail address\/}: tobias\@math.uu.se}\\
\vglue6pt
Stanford University, Stanford, CA, USA\\
{\eightpoint {\it E-mail address\/}: white\@math.stanford.edu}\\
\vglue6pt
Mathematisches Institut, Universit\"at Bonn, 53115 Bonn, Germany\\
{\eightpoint {\it E-mail address\/}:  Daniel\@Wienholtz.de}}

 \vglue12pt
\centerline{\bf Abstract}
\vglue12pt
  This paper proves that classical minimal surfaces 
    of arbitrary topological type with
    total boundary curvature at most $4\pi$ must be smoothly embedded.
    Related results are proved for varifolds and for soap film surfaces.
\bigbreak
 \def\genus{\mathop{\rm genus}}
\def\({\left(}
\def\){\right)}

\def\spt{\mathop{\rm support}}

\def\diverg{\mathop{\rm div}}
\def\ddt{\frac{d}{dt}}
\def\area{\mathop{\rm area}}

\def\diam{\mathop{\rm diam}}

\def\HH{\Cal H}

\def\RR{\Bbb R}
\def\tc{\mathop{\rm Total Curvature}}
\def\Cone{{\textstyle\mathop{\rm Cone}}}
\def\BB{\bold B}
\def\Length{\mathop{\rm Length}}
\def\Len{\mathop{\rm Length}}
\def\kurv{\bold k}
\def\nn{\bold n}
\def\dist{\mathop{\rm dist}}
\def\({\left(}
\def\){\right)}
\def\eps{\varepsilon}
\def\vv{\bold v}
\def\uu{\bold u}
\def\MM{\Cal M}
\def\CC{\bold C}

\def\diverg{\mathop{\rm div}}
\def\pdf#1#2{\frac{\partial #1}{\partial #2}}

\def\GG{\Cal G}

In a celebrated paper \cite{N3} of 1973, Nitsche proved that if $\Gamma$
is an analytic simple closed curve in $\RR^3$ with total 
curvature at most $4\pi$, 
then $\Gamma$ bounds exactly one minimal disk $M$.  Furthermore, 
that disk is smoothly immersed: it has no branch points, either in the
interior or at the boundary.  
His analysis left open the following questions:
\medbreak
\item{(i)} Must $M$ in fact be embedded?
\medbreak \item{(ii)} If $\Gamma$ bounds other minimal surfaces,
  must they also be free of branch points, or even be smoothly embedded?
\medbreak\noindent 
In this paper, we show the answer to both
questions is ``yes,'' even for curves in $\RR^N$.

Regarding \therosteritem{ii}, we give an example of such a $\Gamma$ 
in $\RR^3$ that does indeed bound at least two other minimal surfaces,
namely M\"obius strips.  
We conjecture that any such $\Gamma$ can bound at
most two M\"obius strips, and no surfaces of other topological types.
(See \S 5.)

Before stating our main result, we review some terminology.
The total curvature of a polygonal curve
is the sum of the exterior angles at the vertices.  For an arbitrary
continuous curve, the total curvature is the supremum of the total 
curvatures of inscribed polygonal curves.  This definition, suggested
by Fox, was introduced and analyzed in a paper by Milnor \cite{Mi}.
In case the curve is piecewise smooth, this definition agrees with the 
classical one: the integral
of the norm of the curvature vector with respect to arclength plus
the sum of the exterior angles at the vertices.

Any bounded curve of finite total curvature is rectifiable (i.e., has finite
arclength; see \S10.1).   
If $\gamma$ is an arclength parametrization of such a 
curve, 
then the right and left derivatives $T^+$ and $T_-$
of $\gamma$ exist and are unit vectors
at each interior point (\S10.3,10.4).  
The {\it exterior angle} of the curve at any interior point
is defined to be the angle between $T^+$ and $T_-$.

Finally, if $M$ is a surface in $\RR^N$ and $p$ is a point in $\RR^N$,
the {\it density} of $M$ at $p$ is
$$
   \Theta(M,p) = \lim_{r\to 0} \frac{\area(M\cap\BB(p,r))}{\pi r^2}
$$
provided this limit exists, where $\BB(p,r)$ is the ball of radius $r$
centered at $p$.

We can now state our main theorem:

\nonumproclaim{Theorem {\rm (2.1, 3.2, 4.1)}}
Let $\Gamma$ be a simple closed curve in $\RR^N$ with total 
curvature $\le 4\pi${\rm .}  Let $M$ be a minimal surface with boundary $\Gamma${\rm .}
Then $M$ is embedded up to and including the boundary{\rm ,} with no interior
branch points{\rm .}  
Furthermore{\rm ,} 
at each boundary point $p$ with exterior angle $\theta${\rm ,}
the density of $M$ at $p$ is either $\frac12-\frac{\theta}{2\pi}$ or 
$\frac12 +\frac{\theta}{2\pi}${\rm .}
At a cusp {\rm (}\/i.e.{\rm ,} where $\theta=\pi${\rm ),}
the density of $M$ is $0$ unless $\Gamma$ is contained 
in a plane{\rm .}  
\endproclaim

These density bounds imply, via standard regularity arguments, that for 
piecewise $C^{1,\alpha}$ boundaries, $M$ has a well-defined tangent plane
at every noncusp point, and that in a neighborhood of the point, $M$ is
a graph over the tangent plane.   In particular, if the boundary is smooth,
then the surface is smoothly embedded (and therefore has no branch points.)

Here and throughout the paper, ``simple closed curve'' means
the image of a circle under a continuous one-to-one map.
``Minimal surface''
means classical minimal surface, that is, a continuous and (possibly branched)
conformal harmonic
map of a compact $2$-manifold with boundary into $\RR^N$ such that
the restriction to the boundary is one-to-one.  In Sections 7 and 8, we 
consider more general minimal varieties.
For instance, we prove, roughly speaking, that there are
no singular minimal varieties with total boundary curvature $\le 3\pi$,
and we give a new proof that soap films cannot form on nonclosed
wires with total curvature $\le 2\pi$ \cite{DW}.

The proof that the  interior of $M$ is embedded and unbranched 
(Theorem 2.1) is an immediate
consequence of the following three facts:
\medbreak
\item{(i)} 
The density of a minimal surface $M$ at a point $p\notin\partial M$
is $\le$ the density at $p$ of the cone subtended by $\partial M$.
\smallbreak \item{(ii)}
The density at $p$ of the cone in \therosteritem{i} is at most
$1/(2\pi)$ times the total curvature of $\partial M$.
\smallbreak \item{(iii)} 
The density of $M$ at any branch point or self-intersection
point $p\notin \partial M$ is at least $2$.
\medbreak\noindent 
Fact \therosteritem{i} follows from an 
extension of the familiar monotonicity formula
for minimal surfaces.  We give two proofs of \therosteritem{ii},
one based on the Gauss-Bonnet formula and the other on integral geometry.
Fact \therosteritem{iii} is well-known. 

Analogous facts hold for boundary points.  

After completing this paper, we were informed that fact \therosteritem{i} 
had been
discovered by M. Gromov; see Theorem 8.2.A in \cite{Gr}.
We have retained our proof to keep the paper self-contained and because
(for certain extensions of our main theorem) it is necessary to check that
fact \therosteritem{i} holds under very weak assumptions on boundary behavior.
Also, deducing it from the extended monotonicity formula as we have done
perhaps makes the result more geometrically intuitive (although the proofs
are essentially the same).

The $4\pi$ in our results is in all respects quite sharp.
For example, regarding embeddedness, Almgren and Thurston \cite{AT} and then 
Hubbard  \cite{Hub} (in a simpler way)
showed that
for each $\eps>0$ and $n<\infty$, there is a smooth  
simple closed curve $\Gamma$ in $\RR^3$
(knotted or unknotted as desired)
with total curvature $<4\pi + \eps$ such
that $\Gamma$ bounds within its convex hull 
no embedded surface of genus $\le n$.
Such a curve necessarily bounds nonembedded minimal surfaces of each
genus $\le n$.  
To see this, let $M_g$ be a least-area surface of genus $\le g$
(which exists according to Douglas's theorem: see 
\cite{Sh1}, \cite{DHKW, 11.5}, \cite{TT2}, or \cite{J}).
For $g<n$, $M_g$ necessarily has self-intersections, which implies that its
area may be reduced by cut-and-paste surgery that adds a one-handle.
Thus $\area(M_{g+1})<\area(M_g)$, which implies that the genus
of $M_{g+1}$ is greater than the genus of $M_g$.  Since this holds
for all $g\le n$, in fact we must have $\genus(M_g)=g$ for every
$g\le {n+1}$.

The $4\pi$ is also sharp regarding branch points.  For if 
$F:\BB(0,1)\subset \CC\to \RR^N$ is a minimal surface with a simple
branch point at $0$ 
(i.e., if $\partial F/\partial z:\BB(0,1)\to \CC^N$ 
has order of vanishing $1$ at $z=0$),
then $F$ restricted to a small ball $\BB(0,r)$ will have total boundary 
curvature
only slightly more than $4\pi$.  Of course if $F|\partial \BB(0,r)$ is not
embedded, one can perturb $\BB(0,r)$ slightly to get a domain 
$\Omega\subset\CC$ for which $F|\partial \Omega$ is embedded (except
in the degenerate case when the image of $F$ is contained in a $2$-plane).
Incidentally, $F|\Omega$ is not the only minimal disk bounded by
$F|\partial \Omega$, since the least area disk with this boundary
cannot have branch points (\cite{O1}, \cite{Gu1}, \cite{Alt}).  
Thus this example also shows
that the $4\pi$ in Nitsche's uniqueness theorem is sharp.
Indeed, B\"ohme \cite{Boh} gave examples of curves in $\RR^3$ with
total curvature $4\pi+\eps$ (with $\eps$ arbitrarily small) that bound
many minimal disks, some with branch points.

\bigbreak\centerline{\bf History}\bigbreak

In this section, we survey other theorems asserting that certain
classical minimal surfaces are embedded and/or free of branch points.

The first embeddedness theorem is due to Rad\'o \cite{Ra2, III.10}, 
who proved that
if $\Gamma\subset \RR^N$ has a one-to-one projection onto
the boundary of a convex planar region $R$, then any minimal disk bounded
by $\Gamma$ is the graph of a smooth function over $R$.  In particular, it
is smoothly embedded.
If in addition $N=3$, 
then there is only
one disk and there are no minimal varieties of 
other topological types \cite{Me}.

Concerning branch points, Osserman \cite{O1} proved that an area-minimizing
minimal surface in $\RR^3$ cannot have true interior branch points.
Gulliver \cite{Gu1} and Alt \cite{Alt}
proved that it cannot have false branch points, either.
Gulliver and Lesley \cite{GuL} proved that it 
cannot have boundary branch points
if the boundary curve is analytic.  
It is not known if the Gulliver-Lesley theorem holds for smooth 
(say $C^\infty$) boundaries.   
(See \cite{Wi} for a partial result.)
However, 
Hildebrandt \cite{H} and Heinz-Tomi \cite{HT}
(with subsequent improvements by other authors
\cite{N2}, \cite{K}, \cite{War}, \cite{L})
proved that a minimal surface in $\RR^N$
is as differentiable as its boundary at the boundary itself and that
boundary branch points must be isolated.

Federer \cite{Fed}, refining a theorem of Wirtinger, 
proved that area minimizing surfaces
in $\RR^N$ with $N\ge 4$ can have interior branch points.  
Specifically, he proved that any complex variety in $\CC^N\cong\RR^{2N}$
is absolutely area minimizing: it has less area than any other oriented
surface, regardless of topological type. Thus, for example,
$z\mapsto (z^2,z^3)$ gives a branched, 
area minimizing map from the unit disk
in $\CC \cong \RR^2$ to $\CC^2\cong\RR^4$.
Nevertheless, White \cite{Wh5} showed that the Gulliver-Lesley theorem 
remains true in $\RR^N$.
Here the analyticity is essential: 
there is an example 
of an absolutely area minimizing surface
in $\RR^4$ with a boundary branch point at a $C^\infty$ boundary.
This example is obtained from the example due to 
Gulliver \cite{Gu2} in $\RR^6$ (as stated in \cite{Wh5})
by restricting to the first four coordinates.

Gulliver and Spruck \cite{GuS} proved by a continuity argument 
that if $\Gamma\subset \RR^3$ 
has total curvature $\le 4\pi$
and if it is extreme 
(i.e., it lies on the boundary of a convex set), then
the minimal disk (unique by Nitsche's result) is in fact embedded.
Later Meeks and Yau proved that such a $\Gamma$ does not bound any
other minimal varieties \cite{MY2}.
Gulliver and Spruck 
conjectured that either condition alone (curvature $\le 4\pi$ or 
extreme) suffices for existence of a smooth embedded
minimal disk.
Shortly afterward, 
three sets of authors proved (in very different ways) that
the extremality condition does indeed suffice.
Tomi and Tromba \cite{TT1} proved by a more sophisticated
continuity argument that an extreme 
$\Gamma\subset \RR^3$ bounds a (possibly unstable) smooth embedded minimal
disk.
Almgren and Simon \cite{AS} 
proved that among all smooth embedded disks bounded
by an extreme $\Gamma$, there is one of least area.
Meeks and Yau \cite{MY1} proved the strongest result: 
the Douglas-Rad\'o disk 
(the least area disk among all disks bounded by $\Gamma$) must in fact
be smoothly embedded.
Of course our results show that the $4\pi$ condition also suffices for
embeddedness.

We have already mentioned Nitsche's uniqueness theorem.
The first step in Nitsche's proof was to exclude  
branch points in any minimal disk
with total boundary curvature $\le 4\pi$; this he did 
using
the Gauss-Bonnet formula.
(The Gauss-Bonnet formula for branched surfaces is due
to Sasaki \cite{Sa}, though he made a mistake concerning boundary branch
points. The mistake was pointed out and corrected by Nitsche \cite{N1},
 [N4, \S380].)
Schneider \cite{Sch} also gave a proof 
(using ideas of Rad\'o \cite{Ra2, III})
excluding interior branch points
for curves of total curvature $\le 4\pi$.  
Schneider's proof does not require any smoothness (beyond continuity)
of the boundary curve.
In general, the Nitsche and Schneider arguments 
can be used to bound the total number of
branch points (weighted according to the order of branching) of a minimal 
surface in terms
of the genus and total boundary curvature.  Our approach gives a bound
on the maximum density (and therefore order of branching) at a single point
in terms of the total boundary curvature but independent of genus.

\vglue18pt
\section{Total curvature and densities of cones}

\proclaim{Theorem}
Suppose $\Gamma$ is a closed curve in $\RR^N$ and that 
$p$ is a point not in $\Gamma${\rm .}
Then
$$
   \Length(\Pi_p\Gamma) \le \tc(\Gamma),
$$
where $\Pi_p$ is the radial projection to the unit sphere centered at $p${\rm :} 
$$
\align
 &\Pi_p:\RR^N\setminus \{p\} \to \partial\BB(p,1); \\
 &\Pi_p(x) = p + \frac{x-p}{|x-p|}.
\endalign
$$
Equivalently{\rm ,}
$$
   \Theta(\Cone_p\Gamma,p) \le \frac{\tc{(\Gamma)}}{2\pi},
$$
where
$$
 \Cone_p\Gamma = \{ p + t(x-p): x\in\Gamma, 0\le t\le 1\}.
$$
\endproclaim

\demo{{P}roof {\rm 1 (using Gauss-Bonnet)}}
Consider first the case when $\Gamma$ is smooth.
Since the formula
is not affected by dilations about $p$, 
we may assume that $\Gamma$ lies outside
the unit ball $\BB(p,1)$ centered at $p$.
Let
$$
  A = (\Cone_p\Gamma) \setminus \BB(p,1)
$$
be the annular region between $\Gamma$ and $\Pi_p\Gamma$.
By the Gauss-Bonnet theorem,
$$
  \int_{\partial A} \kurv\cdot \nn\,ds + \int_AK = 2\pi\chi(A),
$$
where $\kurv$ is the curvature vector of the curve $\partial A$ in 
$\RR^N$, $\nn$ is the exterior normal of $A$, $K$ is the scalar
curvature of $A$, and $\chi(A)$ is the Euler characteristic of $A$.
For $A$ we have $K=0$ and $\chi(A)=0$.  Thus
$$
\align
  0 &= \int_{\partial A}\kurv\cdot \nn\,ds    \tag$\ast $\\
\noalign{\vskip4pt}
    &= \int_{\Pi_p\Gamma}\kurv\cdot \nn\,ds   
     +\int_{\Gamma}\kurv\cdot \nn\,ds.       
\endalign
$$  
The first integral in ($\ast$) is equal to the length of $\Pi_p\Gamma$,
and the second is bounded in absolute value by the total curvature of 
$\Gamma$.  This completes the proof when $\Gamma$ is smooth.  
The case of polygonal (i.e., piecewise linear) curves follows 
using approximation
by smooth curves, and the general case then follows by taking the supremum
over inscribed polygonal curves.

(Even when the curve $\Gamma$ is smooth, the annulus $A$ might not be.
However, if we use a point $p'$ in $\RR^{N+1}\supset \RR^N$, the corresponding
annulus $A'$ will be smooth.  Thus the inequality will be true for $p'$
and, letting $p'\to p$, also for $p$.  We thank M. Ghomi for this 
observation.)\hfill\qed
\enddemo 

\demo{{P}roof {\rm 2 (using integral geometry)}}
We may assume $p=0$.  For unit vectors $\vv$, let 
$f_\vv:\Gamma\to \RR$ be the function $f_\vv(x)=\vv\cdot x$.
By Crofton's formula the length of $\Pi_p\Gamma$ is $\pi$ times the
average over 
all~$\vv$ of the number of zeroes of $f_\vv$.
Similarly, Milnor proved that the total curvature of
$\Gamma$ is $\pi$ times the average over all $\vv$  
of the number of local extrema of $f_\vv$ \cite{Mi, 3.1}.

We claim that for each $\vv$, the number of zeroes is $\le$ the 
number of local extrema.  The theorem then follows immediately by averaging
over $\vv$.  To prove the claim, we may suppose $f_\vv$ has finitely many
extrema.  To each zero $q$ of $f_\vv$, associate the first local
extremum $\phi(q)$ that occurs on or after $q$ in the parametrization.
Note that this is an 
injection from the set of zeroes to the set of local extrema. \hfill\qed
\enddemo

\numbereddemo{{R}emark}  
Equality in Theorem 1.1 holds if and only if 
(i) $\Gamma$ lies in a $2$-plane through $p$,
and (ii) within that $2$-plane, $\Gamma$ is locally convex with respect to
$p$.  (Local convexity with respect to $p$ means that
$\Gamma$ is a union of open arcs $U_i$ such that each $U_i$ is in the
boundary of the convex hull of $U_i$ and $p$.)
 
The ``if'' part is easy.
To prove that condition (i) is necessary for equality, suppose 
that $p=0$ and
that $\Gamma$ is not contained in any $2$-plane through $0$.
This implies we can find a small arc $C$ of $\Gamma$ such that
$\Pi_pC$ is contained in an open hemisphere of $\partial \BB(0,1)$, 
but such that $\Pi_pC$ is not contained
in the geodesic joining its endpoints.  For such an arc,
we can find a unit vector $\vv$ such that the minimum of $f_\vv$
on $C$ is positive and strictly less than its value at the two endpoints.
Now (fixing $C$), the set of such $\vv$ is an open set.
For each such $\vv$, note that the number of local extrema of $f_\vv$
is strictly greater than the number of zeroes.  (This is because
the arc of $\Gamma\setminus f_\vv^{-1}(0)$ that contains $C$ has not
just one, but at least three local extrema; the image of the injection
$\phi$ in the proof will only contain at most one of the three.)
Thus when we integrate over all~$\vv$, we get strict inequality.

Necessity of (ii) is proved in a similar way.
\enddemo

\proclaim{Theorem}
Let $M$ be a 
minimal surface in $\RR^N$ with rectifiable boundary $\Gamma${\rm .}
Let $p$ be any point in $\RR^N${\rm .}
Then
$$
   \Theta(M,p) \le \Theta(\Cone_p\Gamma, p)   \tag 1
$$
with strict inequality unless $M = \Cone_p\Gamma${\rm .}
\endproclaim

{\it Remark}.
As mentioned in the introduction, this theorem was discovered by
M. Gromov \cite{Gr, 8.2.A}.

\demo{Proof}
We may assume that $\Theta(\Cone_p\Gamma,p)<\infty$ or, equivalently,
that
$$
   \Len(\Pi_p\Gamma) < \infty.
$$
First consider the case $p\notin \Gamma$.
Let $E=E(\Gamma,p)$ denote the exterior cone over $\Gamma$ with vertex $p$:
$$
   E = \{ p +t(x-p): x\in \Gamma, t\ge 1\}.
$$
Let $M' = M \cup E$ and let
$$
    \Theta(M',p,r) = \frac{\area(M'\cap\BB(p,r))}{\pi r^2}
$$
denote the density ratio of $M'$ in the ball $\BB(p,r)$.  Then
$$
    \Theta(M',p,r) \le \Theta(M',p,R)        \tag 2
$$
for $0<r< R <\infty$.
For $R\le \dist(p,\Gamma)$, this is the standard monotonicity theorem.
For general $R$, this is 
the extended monotonicity theorem 9.1: $\Theta(M',p,r)$ is increasing, 
and is nonconstant unless $M'$
is a cone with vertex $p$.
Letting $r\to 0$ and $R\to \infty$ in \thetag2 gives \thetag1.

In case $p\in \Gamma$, the extended monotonicity formula remains true,
and the proof is exactly as before, except that it is not as obvious
that
$$
    \lim_{r\to 0} \Theta(M',p,r)  = \Theta(M,p).   
$$
To see that this is the case, note that
$$
   \Theta(M',p,r)=\Theta(M,p,r)+\Theta(E,p,r).
$$
Thus we need only show that
$$
   \lim_{r\to 0} \Theta(E,p,r) = 0.
$$
This is intuitively clear if $\Gamma$ is smooth or piecewise smooth.
For the general case, 
note that for $r\le 1$, 
$$
   E\cap \BB(p,r) \subset \Cone_p(\Pi_p(\Gamma\cap \BB(p,r))),
$$
from which it follows that
$$
   \Theta(E,p,r) \le \frac1{2\pi} \Len(\Pi_p(\Gamma\cap \BB(p,r))).  \tag 3
$$
Since $\Pi_p\Gamma$ has finite length, $A\mapsto \Len(\Pi_p\Gamma\vert A)$
defines a finite Borel measure on $G$, where $G$ is a parameter domain
for $\Gamma\setminus\{p\}$.  As $r\to 0$, 
$\Gamma^{-1}(\RR^N\setminus \BB(p,r))$ exhausts $G$ 
and thus $\Len(\Pi_p(\Gamma\cap \BB(p,r)))$ tends to $0$.
Therefore by \thetag3 we get $\Theta(E,p,r)\to 0$ as required.
\enddemo

\section{Interior regularity}

\proclaim{Theorem}
Let $\Gamma$ be a simple closed curve in $\RR^N$ with total curvature
at most $4\pi${\rm .}  Let $M$ be a minimal surface with boundary $\Gamma${\rm .}
Then the interior of $M$ is embedded{\rm ,} and $M$ has no interior branch points{\rm .}
\endproclaim

\demo{Proof}
Since $\Gamma$ has finite total curvature, 
it is also rectifiable (Theorem 10.1),
so we can apply Theorems~1.1 and 1.3.
If $M$ is contained in a cone, then both the mean curvature and
the scalar curvature vanish, so $M$ is contained in a plane.  This case
is trivial, so we assume $M$ is not contained in any cone, which means
we have strict inequality in Theorem 1.3.

Let $p\in M\setminus \Gamma$.  By Theorems 1.1 and 1.3,
$$
  \Theta(M,p) < \frac{\tc{(\Gamma)}}{2\pi} \le \frac{4\pi}{2\pi} = 2. 
$$
Thus $M$ must be unbranched and embedded at $p$, 
since the density of a minimal surface at any interior branch point
or interior self-intersection point is at least~$2$.  
\enddemo

\proclaimtitle{The F\'ary-Milnor Theorem \cite{Fa}, \cite{Mi}}
\proclaim{{C}orollary} 
If $\Gamma$ is a simple closed curve in $\RR^3$ with 
total curvature at most $4\pi${\rm ,} then $\Gamma$ is unknotted{\rm .}
\endproclaim

\demo{Proof}
Let $F:\BB(0,1)\to \RR^3$ be the least area disk (i.e., the Douglas-Rad\'o
solution to the Plateau problem) bounded by $\Gamma$.  By the theorem,
$F$ is an embedding of the interior of $\BB(0,1)$.  In particular,
$r\mapsto F(\partial \BB(0,r))$ 
describes an isotopy of curves for $r\ne 0$.  When $r=1$,
the curve is $\Gamma$, and when~$r$ is\break near $0$, the curve is very nearly 
circular and is therefore unknotted.
\enddemo

 \section{Smooth boundaries}

In this section we exclude boundary branch points for smooth curves
with total curvature at most $4\pi$.  We begin with a theorem about
arbitrary curves.

\proclaim{Theorem}
Let $\Gamma$ be a simple closed curve in $\RR^N$ with finite total
curvature.  Let $p$ be a point in $\Gamma$.  Then
$$
  \Len(\Pi_p\Gamma) \le \tc(\Gamma)-\pi - \theta,
$$
where $\theta$ is the exterior angle to $\Gamma$ at $p${\rm .}
{\rm (}\/Thus $\theta=0${\rm ,} $\theta=\pi${\rm ,} or $0<\theta<\pi$ according to whether
$\Gamma$ has a tangent line{\rm ,} a cusp{\rm ,} or a corner at $p${\rm .)}
\endproclaim

\demo{Proof}
This can be proved by a straightforward modification 
of the integral geometric proof of Theorem 1.1.
Alternatively, we can deduce Theorem~3.1 from Theorem~1.1 as follows.

For small $r>0$, let $a=a_r$ and $b=b_r$ be the intersections of
$\Gamma$ with $\partial \BB(p,r)$.  (There are only two intersection points
by Theorem 10.3.)
Let $\theta(r)$ be the exterior angle of the triangle $apb$ at vertex $p$.  
Let $q$ be a point very close to $p$
such that $p$ is in the triangle $aqb$ but not in the segments $aq$ or
$qb$.
 
Let $\Gamma'=\Gamma'_{r,q}$ be the closed curve obtained by replacing
$\Gamma\cap\BB(p,r)$ with the segments $aq$ and $qb$.
By Theorem 1.1,
$$
  \Len\( \Pi_p\Gamma'_{r,q}  \) \le \tc(\Gamma'_{r,q}).
$$
Note that $\Pi_p\Gamma'$ 
consists of $\Pi_p(\Gamma\setminus \BB(p,r))$ 
together with a geodesic
arc of length $\pi+\theta(r)$, so
$$
  \Len\( \Pi_p(\Gamma\setminus \BB(p,r)) \) + \pi + \theta(r)
        \le \tc(\Gamma'_{r,q}).
$$
Now let $q\to p$:
$$
\align
  \Len\( \Pi_p(\Gamma\setminus \BB(p,r)) \) + \pi + \theta(r)
        &\le \tc(\Gamma'_{r,p}) \\
        &\le \tc(\Gamma).
\endalign
$$
Letting $r\to 0$ gives the result.
\enddemo

\proclaim{Theorem}
Let $\Gamma$ be a smooth simple closed curve in $\RR^N$ with total
curvature at most $4\pi${\rm .}  Let $M$ be a minimal surface with boundary
$\Gamma${\rm .}  Then $M$ is a smooth embedded manifold with boundary{\rm .}
\endproclaim

\demo{Proof}
We already know from Theorem 2.1 
that the interior of $M$ is smooth and embedded.
Let $p\in \Gamma$.  As in the proof of Theorem 2.1, 
we may assume we have strict inequality in Theorem 1.3.  
Then by Theorems 1.3 and 3.1, \pagebreak
$$
\Theta(M,p) 
< \Theta(\Cone_p\Gamma,p)  
\le  \frac{\tc(\Gamma)}{2\pi} - \frac12
\le \frac{4\pi}{2\pi} - \frac12 = \frac32.     
$$
But any boundary branch point has density at least $3/2$.
Likewise, if $M$ were not embedded at $p$, then the density would be
at least $3/2$.  
\enddemo

\section{Boundaries with corners}

\proclaim{Theorem}  Let $\Gamma$ be a simple closed
curve in $\RR^N$ with total curvature at most $4\pi${\rm .}  
Let $M$ be a minimal surface with boundary $\Gamma${\rm .}
\medbreak
\item{\rm (i)} If $p$ is a point of $\Gamma$ with exterior angle $\theta${\rm ,}
then the density $\Theta(M,p)$ is either 
$\frac12+\frac{\theta}{2\pi}$ or $\frac12 -\frac{\theta}{2\pi}${\rm .}
\smallbreak \item{\rm (ii)} If $p$ is a cusp{\rm ,} then the density of $M$ at $p$ is $0$ unless
$\Gamma$ lies in a plane{\rm .} \smallbreak
\item{\rm (iii)} $M$ is embedded up to and including the boundary{\rm .}
\endproclaim 

\demo{Proof}
Let $T$ be a tangent cone to $M$ at $p$.  Then the boundary
of $T$ consists of two rays with an interior angle of $\pi-\theta$.
The intersection $S$ of $T$ with the 
unit sphere is a collection of one or more
geodesic curves (see  Note~4.2 below).  
One curve must be a geodesic arc joining
two points that are $\pi-\theta$ apart (in geodesic distance).
Thus this arc must have length $\pi-\theta$ or $2\pi-(\pi-\theta)=\pi+\theta$.
The other arcs must be great circles.  Thus the density of $M$ at $p$
must be
$$
   \Theta(M,p) =  \frac1{2\pi}(\pi \pm \theta + 2\pi k),  \tag 1
$$
where $k\ge 0$ is the number of great circles.

Now by Theorems 1.3 and 3.1,
$$
  \Theta(M,p) < \frac1{2\pi}(\tc(\Gamma)-\pi-\theta)
              \le \frac1{2\pi}(3\pi - \theta).             \tag 2
$$
(The strict inequality comes from assuming that $M$ is not contained
in a cone; otherwise all the conclusions are trivially true.)
Combining \thetag1 and \thetag2 we see that
$$
   \pi \pm \theta + 2\pi k \le 3\pi - \theta
$$
or
$$
   (\theta \pm \theta) + 2\pi k < 2\pi.    \tag 3
$$
Notice this forces $k=0$.  If $p$ is a cusp (so $\theta=\pi$),
then the $\pm$ in \thetag3 (and therefore in \thetag1) 
must be $-$, which means that (in this case) 
the density $\Theta(M,p)$ is $0$.

Note that if $p\in \Gamma$ coincided with an interior point of $M$,
then $S$ would indeed contain at least one great circle, which we
have just seen it does not.  This proves the embeddedness.
\enddemo

\numbereddemo{Note} 
In the proof of Theorem 4.1 (and later in 7.1 and 7.2), 
we use the following well-known fact: 
any tangent cone to a $2$-dimensional 
minimal variety (e.g. stationary integral varifold) 
intersects the unit sphere in a finite collection of geodesic arcs.  
This fact can be proved as follows.
The intersection of the tangent cone with the unit sphere is
a $1$-dimensional varifold that is, by a straightforward calculation,
stationary in the sphere.
By \cite{AA}, such $1$-dimensional stationary
varifolds consist of geodesic arcs.
\enddemo
 
\section{Nondisk type surfaces}

In this section we prove that there are simple closed curves
in $\RR^3$ with total curvature $<4\pi$ that bound minimal M\"obius
strips.

Consider two convex polygons in the halfplane $\{(x,y,0):y\ge 0\}$
such that each polygon has exactly one vertex, namely the origin, located
on the $x$-axis.  (Think of the two polygons being the same regular $n$-gon.)
Give both polygons the counterclockwise orientation.
Now rotate one of the two polygons about the $x$-axis through a small 
positive angle and the other through a small negative angle.
The two polygons will still have the origin as a common vertex.

Consider the closed connected curve $\Gamma$ that, 
starting at the origin, traces
out one polygon and then the other.
We claim that $\Gamma$ has total curvature $<4\pi$.
To see why, note that
for a generic unit vector $\vv$, the function $f_\vv(x)=\vv\cdot x$
will not be constant on any segment of $\Gamma$.
For such a $\vv$, $f_\vv$ can have at most four  local 
extrema (two on each polygon).  
However, the set of such $\vv$ for which $f_\vv$ has only two local
extrema is open and nonempty since it contains vectors arbitrarily close
to $(0,0,1)$.  Thus the average number of local extrema is less than $4$,
which implies \cite{Mi, 3.1}
that the total curvature is less than $4\pi$.

Of course $\Gamma$ has a self-intersection, 
but we can move one vertex slightly to make it embedded.
We can also round the corners to make it smooth or, 
by further perturbation,
even analytic.   

Note we can construct such a $\Gamma$ that, after a suitable
translation, is arbitrarily close to 
a curve that traverses the unit circle in the $xy$ plane twice.  
We can also make $\Gamma$ lie outside the unit cylinder $\BB(0,1)\times\RR$.
Fix such a curve $\Gamma$.

Now a disk (or any other orientable surface)
bounded by $\Gamma$ has area $\ge 2\pi$, since its projection to the
$xy$-plane must cover the unit disk twice algebraically.
But clearly $\Gamma$ bounds a M\"obius strip with much smaller area.
Thus by the unoriented version of the Douglas theorem 
(see \cite{Sh1} or \cite{Be}), $\Gamma$ bounds a minimal M\"obius
strip $M$.

If we perturb $\Gamma$ slightly, then it will
not bound any minimal surfaces with nonzero nullity (by Theorem 11.1(i)).  
Hence $M$ will be strictly
stable.  It then follows from degree theory (\S11.1(ii))
that $\Gamma$ must also bound an unstable minimal M\"obius strip.

By analogy with certain results and conjectures \cite{MW}
about minimal surfaces
bounded by pairs of convex curves, we make the following conjecture.

\nonumproclaim{{C}onjecture} Let $\Gamma$ be a smooth simple closed curve
in $\RR^N$ with total curvature $\le 4\pi${\rm .}  Then in addition to a unique
minimal disk{\rm ,} $\Gamma$ bounds either\/{\rm :}\/
\medbreak
\item{\rm (i)} no other minimal surfaces{\rm ,} or
\smallbreak \item{\rm (ii)} exactly one minimal M{\rm \"o{\it }}bius strip and no
  other minimal surfaces{\rm ,} or
\smallbreak \item{\rm (iii)} exactly two minimal M{\rm \"{\it o}}bius strips and no other minimal surfaces{\rm .}
\medbreak\noindent 
In case {\rm \therosteritem{ii},} the strip has index $0$ and nullity $1${\rm .}  In
case {\rm \therosteritem{iii},} both strips have nullity $0${\rm ,} one has index $1$ and
the other has index $0${\rm .}
\endproclaim

We emphasize that this conjecture refers to 
classical minimal surfaces bounded by $\Gamma$:
in Section 7 we show that $\Gamma$ can also bound  
soap-film-like minimal varieties.

It would be interesting to know the least possible total boundary curvature
of a minimal M\"obius strip, or even a good lower bound.  One could
ask the same question about higher genus nonorientable surfaces, although,
if the above conjecture is correct, the answer would be $4\pi$ by the
Almgren-Thurston-Hubbard examples.

\section{Disconnected boundaries}

Suppose $M$ is a minimal surface with more than one boundary component,
and suppose that the total curvature of the boundary (i.e., the sum
of the total curvatures of the components) is at most $4\pi$.
By the Fenchel-Borsuk theorem (\cite{Fen}, \cite{Bor}, \cite{Mi}, \cite{Fa}), 
each component has total curvature $\ge 2\pi$
and indeed has total curvature $>2\pi$ unless it is a plane convex curve.
Thus $\partial M$ must consist of exactly two components $\Gamma_1$ and
$\Gamma_2$, each of which is a plane convex curve.  

If $M$ is a cone, then its scalar curvature and mean curvature
both vanish, so it is locally planar.  That is, it is the union of the two
planar regions $R_1$ and $R_2$ bounded by $\Gamma_1$ and $\Gamma_2$.
In this case, the two regions must intersect since the vertex of the cone
belongs to both regions.
If $M$ is not a cone, then all the conclusions of Theorems 2.1, 3.2, and 4.1
hold, with exactly the same proofs.  (Connectivity of $\partial M$ was
only used to conclude that $M$ was not a cone.) \pagebreak

Meeks and White \cite{MW} studied the case when $\Gamma\subset \RR^3$
and $R_1$ and $R_2$ are disjoint.
They proved, for example, that such a
$\Gamma$ bounds
at most two minimal surfaces of annular type, and that only one of the annuli
can be stable.  
They also conjectured that such a $\Gamma$ cannot bound minimal surfaces
of higher genus.
If, on the other hand, $R_1$ and $R_2$ are not disjoint but rather intersect
along a segment, then $\Gamma$ will bound other surfaces.  In particular,
there are examples \cite{MW, 3} that bound more than one stable annulus.
Futhermore, 
by the heuristic principle that every minimal surface with a curve
of self-intersection
should be the limit (as a set) of a sequence of embedded minimal surfaces,
we expect that such a $\Gamma$ should bound embedded minimal surfaces of
arbitrarily large genus.  (Apply the heuristic principle to $R_1\cup R_2$.)

\section{Nonclassical minimal surfaces}

As we have seen, a smooth simple closed curve of finite total curvature
$\le 4\pi$ does not bound any classical minimal surface with branch points
or with self-intersections.
However, it may bound other minimal varieties, such as stationary integral
varifolds and flat chains modulo $k$, that do have singularities.

For example, consider a curve $\Gamma_0$ in $\RR^3$ 
that traverses the unit circle in the $xy$ plane twice
and therefore has total curvature $4\pi$.  
As shown in Section 5, one can modify the curve slightly
to get an embedded curve $\Gamma$ with total curvature 
less than $4\pi$ that lies just outside the cylinder $\BB(0,1)\times\RR$.

Let $T$ be a mod $3$ area minimizing surface (i.e., flat chain)
with boundary $\Gamma$.
If $T$ were smooth, it would be orientable (since $1\not\equiv -1$ mod $3$),
and so its area would be $\ge 2\pi$ since its projection to the $xy$
plane would have to cover the unit disk twice algebraically.

However, $T$ has area at most
slightly greater than $\pi$.  
To see this, consider the surface $T'$ consisting of a ribbon
homotoping $\Gamma$ to $\Gamma_0$, together with the disk $\BB(0,1)$ with 
multiplicity $-1$.  Then the boundary of $T'$ is congruent to $\Gamma$ mod 
$3$, and the mass is only slightly more than $\pi$.
Thus $T$ cannot be smooth.  

Now consider a rectifiable curve $\Gamma$ in $\RR^N$.
Recall that a compactly supported rectifiable varifold $M$ is 
said to be stationary in $\RR^N\setminus\Gamma$
if
$$
   \left(\ddt\right)_{t=0} \area(\phi(t,M)) = 0
$$
for all smooth $\phi:\RR\times\RR^N\to\RR^N$ such that $\phi(0,x)\equiv x$
and $\phi(t,\xi)\equiv \xi$ for $\xi$ in a neighborhood of 
$\Gamma$. 
Let us call $M$ {\it strongly stationary} with respect to $\Gamma$ provided
$$
   \left(\ddt\right)_{t=0}
   \left( \area(\phi(t,M)) + \area(\phi([0,t]\times\Gamma))  \right)
        \ge 0
$$
whenever $\phi:\RR\times\RR^N\to \RR^N$ is a smooth map 
with $\phi(0,x)\equiv x$.  Intuitively, strong stationarity means that
the area of $M$ cannot be reduced to first order by any smooth deformation,
even moving the boundary, 
provided we include the ribbon swept out by $\Gamma$
as part of the deformed surface.

Strong stationarity may also be stated in any of the following three
equivalent ways:
\medbreak
\item{(i)} For any smooth map $\phi:\RR\times\RR^N\to \RR^N$ with
$\phi(0,x)\equiv x$, we have
$$
    \(\ddt\)_{t=0} \area(\phi(t,M)) \le \int_\Gamma|X^\perp|\,ds,
$$
where $X(x)=\(\ddt\)_{t=0}\phi(t,x)$ and $X^\perp$ is the component
of $X$ normal to $\Gamma$.
\smallbreak \item{(ii)} $\int_M\diverg_M X \le \int_\Gamma |X^\perp|\,ds$ 
for smooth vectorfields
    $X$ on $\RR^N$.
\smallbreak \item{(iii)} 
There is an $\HH^1$-measurable normal vectorfield $\nu$ on $\Gamma$
with $\sup|\nu|\le 1$
such that $\int_M\diverg_M X = \int_\Gamma X\cdot \nu\,ds$
for all smooth vectorfields $X$ on $\RR^N$.
\medbreak\noindent
The integrals over $\Gamma$ are with respect to arclength or, 
equivalently, one-di\-mensional Hausdorff measure.
The equivalence of strong stationarity and \therosteritem{i} is clear.
Statements \therosteritem{i} and \therosteritem{ii} are equivalent by
the first variation formula, and \therosteritem{ii} and \therosteritem{iii}
are
equivalent by the Riesz representation theorem.

We emphasize (as evidenced by \therosteritem{ii} and \therosteritem{iii})
that strong stationarity is really stationarity in $\RR^N\setminus\Gamma$
plus a boundary condition.  
 Indeed, for $C^{1,1}$ curves $\Gamma$, 
the boundary condition is that the first variation
boundary measure is $\le$ arclength measure on $\Gamma$.  (This would
be true for general rectifiable $\Gamma$ 
if the varifold-theoretic unit normal
at the boundary were known to be perpendicular to the boundary almost 
everywhere. For the $C^{1,1}$ case, see \cite{All2}.)

Thinking physically, one could say that $M$ is strongly stationary
if it is stationary, if the total force per unit length that the surface
exerts on $\Gamma$ is everywhere $\le 1$, and if this force is perpendicular
to the boundary.

To better understand this definition, note that
any classical minimal surface with
smooth boundary $\Gamma$ is (as a varifold) strongly stationary with respect
to $\Gamma$.  The union of two such surfaces with the same boundary
will not in general define a strongly stationary varifold.  
Indeed, the resulting varifold is strongly stationary if and only
if the angle between the the two surfaces is everywhere $\ge 120$ degrees
along the boundary.
(To see this, use the first variation formula.)
This corresponds to what one observes in soap films: two sheets
can meet along a portion of wire, but never with an angle of less than
$120$ degrees.

\proclaim{Theorem}
Let $\Gamma$ be a simple closed curve in $\RR^N$ with total curvature
at most $3\pi${\rm .}  Suppose $M$ is a compactly supported rectifiable varifold
that is strongly stationary with respect to $\Gamma${\rm .}
Suppose also that the density of $M$ is $\ge 1$ at every point
in $(\spt M)\setminus \Gamma${\rm .}
Then $M$ is smooth in the interior{\rm ,} and also at the boundary if the 
boundary is smooth{\rm .}
\endproclaim

\demo{Proof}
As in the proof of Theorem 2.1,
we can assume that $M$ is not contained in a cone.
Consider first a point $p\notin \Gamma$.  
Then
$$
   \Theta(M,p) < \Theta(\Cone_p\Gamma,p) 
      \le \frac1{2\pi}\tc(\Gamma) \le \frac32.
$$
(The first inequality is the varifold version of Theorem 1.3:
the proof is identical, except that one uses 
the varifold version (\S9.2) of the extended monotonicity formula.)
But this density bound implies that $p$ is regular as follows:

Let $C$ be a tangent cone at $p$.  Then $\Theta(C,0)<3/2$.
Thus $\Theta(C,x)<3/2$ for all $x$, since in any minimal cone,
the highest density occurs at the vertex.  
(This is because the density function
$\Theta(C,\cdot)$ is upper-semicontinuous   [Si, \S17.8]
and constant along radial lines.)
Now the intersection of $C$ with the unit sphere is a collection
of geodesic arcs (\S4.2), which means that the cone $C$ is polyhedral.
At most two faces of $C$ can meet along an edge, since otherwise
the density at points in the edge would be $\ge 3/2$. 
This means $C\cap \partial \BB$ is a union of great circles. Since
the density is $<3/2$, there is only one great circle and it has multiplicity
$1$. By Allard's regularity theorem (\cite{All1} or  [Si, \S24.2]),
this means that $M$ is regular at $p$.

The boundary case is similar: one shows that the density at any boundary
point is $<1$, whereas the minimum possible density at a boundary singularity
is $1$.
\enddemo

This theorem becomes false if $3\pi$ is replaced by $4\pi$: the varifold
associated to the mod $3$ flat chain $T$ discussed at the beginning of this
section is a counterexample.

The surface $M$ in Theorem 7.1 is not assumed to 
minimize area in any sense.
For surfaces in $\RR^3$ that minimize area modulo $k$,
less control on total boundary curvature is needed to ensure
regularity:

\proclaim{Theorem} Let $\Gamma$ be an embedded closed curve
{\rm (}\/not necessarily connected\/{\rm )} in $\RR^3$ with 
total curvature $\le k\pi${\rm .}
Let $M$ be a mass\/{\rm -}\/minimizing surface {\rm mod} $k$ with boundary $\Gamma${\rm .}
Then $M$ is smooth and embedded in the interior{\rm ,} and also at the boundary
if the boundary is smooth{\rm .}
\endproclaim

\demo{Proof} 
Consider a tangent cone $C$ to $M$ at a point $p\notin\Gamma$.
As in the proof of Theorem 7.1, $C$ must be a
polyhedral complex.  If the complex has edges, then $k$ (or a multiple of 
$k$) sheets meet
along such an edge which (as in the beginning of this section
for the case $k=3$) means that $C$
has density $\ge k/2$.  Thus if $C$ has density $<k/2$, it has no edges and
is therefore a plane
with some multiplicity $<k/2$.\pagebreak

Now, as in the proof of Theorem 7.1, one deduces that the density at each
interior point $p$ of $M$ is $<k/2$.  Let $C$ be a tangent cone
to $M$ at $p$.  Then, as we have just seen, $C$ is a plane of multiplicty
$<k/2$.  But that implies that $p$ is a regular point of $M$ 
\cite{Wh1, Theorem 3.1}.

Thus the interior of $M$ is a smooth minimal surface.
Since the density is everywhere $<k/2$, $M$ is orientable.
Consequently $M$ is mass-minimizing as an integral current.  Thus
it is smooth at the boundary wherever the boundary is smooth
\cite{HS}.
\enddemo

Does this theorem also hold in $\RR^N$ for $N>3$?
The proof of the main theorem of \cite{Wh1} 
is valid only for hypersurfaces.  However, that theorem is probably
true for arbitrary codimensions.  
If so, then Theorem 7.2 does indeed hold
for surfaces in $\RR^N$.

Probably neither of these last two theorems (7.1 or 7.2) is sharp.
However, as $k\to\infty$ the ratio of the bound in Theorem 7.2 
to the optimal bound
tends to $1$.  
For consider an embedded curve that is a slight perturbation
of a circle traversed $(k+1)/2$ times if $k$ is odd or $(k+2)/2$ 
times if $k$ is even.  Then the total curvature is slightly less than
$(k+1)\pi$ or $(k+2)\pi$, respectively, 
and the mass minimizing surface mod $k$
will have singularities.  The proof is just like the proof for the
special case $k=3$ given at the beginning of this section.

Morgan \cite{Mo1} proved an interesting theorem 
related to Theorem 7.2: 
if an $m$-dimensional mass
minimizing hypersurface mod $p$ has a smooth, extreme boundary
with at most $p/2$ connected components, then the hypersurface
is smooth except for a set of Hausdorff dimension $\le m-7$.

\section{Nonclosed curves}

\proclaim{Theorem} Let $\Gamma$ be a smooth embedded connected
nonclosed curve in $\RR^N$ with total curvature $\le 2\pi${\rm .}
Then no nontrivial varifold $M$ with bounded support and with
density $\ge 1$ throughout $(\spt M)\setminus \Gamma$ is strongly
stationary with respect to $\Gamma${\rm .}
\endproclaim

This is sharp because for every $\eps>0$, there are examples (see \cite{DW})
with total
curvature $\le 2\pi +\eps$.  Readers can construct such examples 
experimentally by tying a simple overhand knot in a piece of wire and
dipping it in soapy water.  

\demo{Proof}
Suppose $M$ were such a varifold.  
We may assume that $\Gamma$ is compact; otherwise replace $\Gamma$
by the smallest connected subset of $\Gamma$ that contains 
$\Gamma\cap \spt M$.
We claim that
$$
 \Theta(\Cone_p\Gamma,p) \le \frac{\tc(\Gamma) + \pi}{2\pi}  \tag$\ast$
$$
for any nonclosed curve $\Gamma$ and point $p\notin\Gamma$.
Using this (and the analogous inequality at the boundary), 
one shows exactly as in the proof of Theorem~7.1 that $M$ is a smooth embedded manifold
with boundary.  But then the boundary must be a closed curve.
Since $\Gamma$ contains no closed curve, this means that the boundary
of $M$ is empty.  But there are no compact, boundaryless minimal surfaces
in $\RR^N$, a contradiction.

To see ($\ast$), attach long line segments to the ends of $\Gamma$ to get
a new curve with endpoints in a large sphere $\BB(p,R)$ centered at $p$
and with the same total curvature.
Now join the ends by a geodesic arc in $\BB(p,R)$
to get a closed curve~$\Gamma'$.
The total curvature of $\Gamma'$ will exceed that of $\Gamma$ by approximately
$\pi$ (for the two new corners) plus the angle $\alpha$
subtended by the geodesic arc.
The length of $\Pi_p\Gamma'$ will exceed that of $\Pi_p\Gamma$ by 
at least $\alpha$.
Hence applying Theorem~1.1  (or Theorem 3.1 if $p\in\Gamma'$) 
to $\Gamma'$ and then letting 
$R\to \infty$\break gives ($\ast$).
\enddemo

Drachman and White \cite{DW} proved (when $N=3$) the stronger result
that there is no compactly supported rectifiable varifold that
is stationary in $\RR^N\setminus \Gamma$.
%
%
%

\section{The extended monotonicity theorem}

\proclaim{Theorem} 
Let $M$ be a minimal surface in $\RR^N$ such that $\Gamma=\partial M$
has finite length{\rm ,} and let $p$ be a point in $\RR^N${\rm .}  
Let $M'$ be the union of $M$ and the exterior cone $E$ over $\Gamma$
with vertex $p${\rm .}  Then
$$
  \Theta(M',p,r) := \frac{\area(M'\cap\BB(p,r))}{\pi r^2}  \tag*
$$
is an increasing function of $r${\rm ,} and is not constant unless $M'$
is a cone{\rm .}
\endproclaim

\demo{Note} 
If $\gamma:S^1\to \RR^N$ is a Lipschitz parametrization of $\partial M$
(e.g., para\-metrization by arclength), then $E$ is given by the mapping
$$
\align
   &e: S^1\times [1,\infty) \to \RR^N, \\
   &e(s,\lambda) = \lambda \gamma(s).
\endalign
$$
Of course the area referred to in ($\ast$) is mapping area (or, equivalently,
area of the image counting multipicities).
\enddemo

In the case of area-minimizing surfaces $M$, the extended monotonicity 
theorem appeared in \cite{Wh4}.

\demo{Proof} The proof of the monotonicity formula in 
[Si, \S17] for manifolds without boundary also works for manifolds with 
boundary, provided a boundary term is added. More precisely, 
suppose $\Sigma \subset \RR^N$ is a smooth submanifold with boundary.
Let $H$ be its mean curvature vector, $\nn_\Sigma(x)$ be the 
outward pointing unit normal to $\partial \Sigma$ at $x$ 
(that is, $\nn_\Sigma(x)$ is tangent to $\Sigma$, normal to 
$\partial \Sigma$, and points out of $\Sigma$), and let
$\rho(x)=|x-p|$.  Then
$$ 
\align 
  \frac{d}{dr}(\Theta(\Sigma, p,r)) 
  = 
  & \frac{d}{dr}\int_{\Sigma\cap \BB(p,r)} \frac{|D^\perp \rho|^2} 
    {\rho^2}\, dA  \tag 1 
  \\ 
  &+ 
  r^{-3} \int_{\Sigma\cap \BB(p,r)} (x-p)\cdot H \, dA 
  -
  r^{-3} \int_{\partial \Sigma \cap \BB(p,r)} (x-p)\cdot \nn_\Sigma 
    \, ds. 
 \endalign 
$$ 
On the other hand, straightforward calculation shows 
that if $\Sigma$ is the exterior cone with vertex $p$
over a piecewise smooth curve 
$\partial\Sigma$, then 
$$ 
  \frac{d}{dr} (\Theta(\Sigma,p,r)) 
  = 
  - r^{-3} \int_{\partial \Sigma \cap \BB(p,r)} (x-p)\cdot \nn_\Sigma 
    \, ds. 
  \tag 2 
$$ 

Now, consider first the case that $\Gamma$ is smooth or piecewise smooth. Then 
we can apply formulas \thetag 1 and \thetag 2 to $M$ and $E$, respectively. 
Note the term involving $H$ vanishes for $M$ since it is minimal. 
Thus adding the formulas gives
$$ 
\align 
  \frac{d}{dr}(\Theta(M',p,r)) 
  = 
  & \frac{d}{dr} \int_{M\cap \BB(p,r)} \frac{|D^\perp \rho|^2}{\rho^2} 
    \, dA   \tag 3 
  \\ 
  & - r^{-3} \int_{\Gamma\cap \BB(p,r)} (x-p)\cdot(\nn_M + \nn_E)\,ds . 
\endalign 
$$ 
Let $T(x)$ be a unit tangent vector to $\partial M$ at $x$.  
Consider the quantity
$$
    (x-p)\cdot \nn 
$$
as a funtion of unit vectors $\nn$ perpendicular to $T(x)$.
Note that the minimum and maximum occur at $\nn=\nn_E$ and $\nn=-\nn_E$,
respectively.  Hence 
$$
  (x-p)\cdot(\nn_M+\nn_E) \le (x-p)\cdot(-\nn_E+\nn_E) = 0
$$
thus \thetag3 becomes:
$$
  \frac{d}{dr}(\Theta(M',p,r))
  \ge
  \frac{d}{dr}\int_{M\cap \BB(p,r)}\frac{|D^\perp \rho|^2}{\rho^2}\,dA
$$ 
or
$$
  \Theta(M',p,R)-\Theta(M',p,r)
 \ge
  \int_{M\cap (\BB(p,R)\setminus\BB(p,r))}
         \frac{|D^\perp \rho|^2}{\rho^2}\,dA.
\tag 4
$$

This proves the result if $\partial M$ is smooth or piecewise smooth.

Now suppose that $\partial M$ is merely rectifiable.  
Then (see the note below)
we can exhaust $M$ by subdomains $M_i$ such that $\partial M_i$
is smooth, $\partial M_i$ converges to $\partial M$ uniformly, and
the length of $\partial M_i$ converges to the length of $\partial M$.

It follows that $M_i$ and $E_i$ converge to $M$ and $E$ as measures
and therefore that $\Theta(M_i',p,r)$ converges to $\Theta(M',p,r)$
for almost every $r$.  
Fix $r$ and $R$ with $r<R$ for which the convergence holds.
By \thetag4 applied to $M_i'$,
$$
\align
  \Theta(M_i',p,R)-\Theta(M_i',p,r)
 &\ge
  \int_{M_i\cap (\BB(p,R)\setminus\BB(p,r))}
         \frac{|D^\perp \rho|^2}{\rho^2}\,dA    \\
 &\ge
  \int_{M_j\cap (\BB(p,R)\setminus\BB(p,r))}
         \frac{|D^\perp \rho|^2}{\rho^2}\,dA    
\endalign
$$
for $j<i$ since then $M_j\subset M_i$.
Letting first $i$ and then $j$ tend to infinity gives \thetag4 for $M'$
itself.
\enddemo

\demo{Note}
If $M$ is a disk, then the fact that 
$\partial M$ can be approximated in $M$ by a smooth curve
of nearly the same arclength is due to Rad\'o \cite{Ra1}.
(Rad\'o's proof is also given in  [N4, \S316].)
For general $M$, we indicate a quick proof similar to Rad\'o's.  
Let $R$ be an annular region
in $M$, one component of which is a connected component
of $\partial M$ and the other a smooth curve.  Then $R$ can be
parametrized conformally by a map
$$
   F: (x,y)\in A=S^1\times [0,a] \to \RR^N.
$$
Of course since $M$ is minimal, $F$ will also be harmonic.
We claim that the arclength $L(y)$ 
of $F(\cdot,y)$
is a convex function of $y$. 
Now $\pdf{F}x$ is harmonic, so $|\pdf{F}x|$ is subharmonic.  (Of course
at the boundary of $A$, this has to be understood as a distribution or
measure.)
Note that $(x,y)\mapsto L(y)$ is the average of the subharmonic functions
$|\pdf{F}{x}(x+\theta,y)|$ (averaged over $\theta$) and hence is subharmonic.
But $L(y)$ depends only on $y$, so subharmonicity is equivalent to convexity.

Distribution derivatives at the boundary may be avoided by replacing
$|\pdf{F}x|$ by the subharmonic function 
$\sum_{j=0}^n|F(w_n^{j+1}x,y)-F(w_n^jx,y)|$ where $w_n=\exp(2\pi i/n)$, 
averaging as above, and then letting $n\to\infty$.

Convexity of the function $L(y)$ was first proved, using a somewhat more
complicated argument, by Osserman and Schiffer \cite{OS}.  They also proved
the sharper result 
that $L''(y) \ge L(y)$, with equality precisely for catenoids.
\enddemo

 9.2. {\it Varifolds}.
The extended monotonicity theorem is also true for $m$-dimensional
rectifiable
varifolds $V$ that are strongly stationary (see \S 7) with respect
to a closed rectifiable set $\Gamma$.
Of course if $m\ne 2$, then the expression $\pi r^2$ in 9.1 ($\ast$) 
should be replaced
by $\omega_mr^m$ where $\omega_m$ is the volume of the unit ball
in $\RR^m$.
The proof is exactly as above.  The exterior cone should be interpreted
as the varifold image of the multiplicity $1$ varifold
associated with $\Gamma\times [1,\infty)$ under the map $e$ defined 
at the beginning of this section.\pagebreak

\demo{{\rm 9.3.} Mass minimizing flat chains}
The extended monotonicity theorem is also true for $m$-dimensional
mass minimizing integral currents or, more generally, for mass minimizing
flat chains with coefficients in any normed group.
The exterior cone should be interpreted as the 
flat chain image $e_\#(\partial M \times [1,\infty))$, where
$e$ is as in the note at the beginning of this section.
As in \S9.2, 
if $m\ne 2$, then the expression $\pi r^2$ in 9.1 ($\ast$) 
should be replaced
by $\omega_mr^m$, where $\omega_m$ is 
the volume of the unit ball in $\RR^m$.
The proof is exactly like the proof
by cone comparison (cf. \cite{Mo2, 9.3})
of the usual monotonicity formula. \enddemo

\section{Two basic properties 
               of curves with finite total curvature}

{\it Rectifiability}.

\proclaim{Theorem}
If $\Gamma$ is a compact connected curve in $\RR^N$ with
finite total curvature{\rm ,} then it has finite length{\rm .}
\endproclaim

\demo{Proof} Let $\gamma(t)$ be a parametrization of $\Gamma$.
We may assume that $\gamma$ is closed (otherwise close it up
with a straight line segment).  If $\uu$ is a unit vector, then the total
variation of $t\mapsto \gamma(t)\cdot \uu$ is at most the diameter
of $\Gamma$
times the number of local extrema of $\gamma(t)\cdot \uu$.
Averaging both sides of this inequality (over all unit vectors $\uu$)
gives
$$
  c_N\Len(\Gamma) \le  (\diam\Gamma)\cdot \tc(\Gamma),
$$ 
where $c_N>0$ depends only on the dimension $N$.
(Here, as in the proof of Theorem 1.1, 
we are using the integral geometric formulas
of Crofton and Milnor \cite{Mi, 3.1}.)
\enddemo

{\it Existence of strong one-sided tangents}.

\proclaim{Lemma} 
Suppose $\gamma:[a,b] \to \RR^N$ is an injectively parametrized
curve of finite total curvature{\rm .}
For $\xi<\eta$ in $[a,b]${\rm ,} let
$$
   T_{\xi\eta} = \frac{\gamma(\eta)-\gamma(\xi)}{|\gamma(\eta)-\gamma(\xi)|}.
$$
If $a< x \le y < b${\rm ,} then
the angle $\angle(T_{ax},T_{yb})$
between $T_{ax}$ and $T_{yb}$ is less than or equal to the total
curvature of $\gamma|(a,b)${\rm .}
\endproclaim

\demo{Proof}
$$
\align
\tc(\gamma|(a,b))
&\ge
\angle(T_{ax},T_{xy})+\angle(T_{xy},T_{yb}) \\
&\ge
\angle(T_{ax},T_{yb})
\endalign
$$
by the triangle inequality for geodesic distances in the unit sphere.
\enddemo

\proclaim{Theorem}
Suppose $\gamma:[A,B] \to \RR^N$ is an injectively parametrized
curve of finite total curvature{\rm .}
Then the strong one\/{\rm -}\/sided tangents
$$
T^+(a) = \lim\Sb a\le x<y \\ y\to a\endSb T_{xy}
\quad
\text{and} 
\quad
T_-(b) = \lim\Sb x<y\le b \\ x\to b\endSb T_{xy}
\tag 1
$$
exist for every $a\in [A,B)$ and $b\in (A,B]${\rm .}  
Furthermore{\rm ,} $T^+(x)$ and $T_-(x)$ both approach $T^+(a)$
as $x$ approaches $a$ from the right {\rm (}\/i.e{\rm ,} with $x>a${\rm ),} and they both
approach $T_-(b)$ as $x$ approaches $b$ from the left {\rm (}\/i.e.{\rm ,} with $x<b${\rm ).}
\endproclaim

\demo{{R}emark} After completing this paper, we discovered that this lemma
appears as lemma 2.13 in \cite{Ku}.  We include the proof below to keep
the paper self-contained and because it seems to us more direct than the
proof in \cite{Ku}.
\enddemo

\demo{Proof}
Let $T$ be a subsequential limit of $T_{ax}$ as $x\to a$ with $x>a$.
Thus by lemma 10.2 applied to $\gamma\vert [a,y]$,
$$
  \angle(T,T_{xy}) \le \tc(\gamma|(a,y))       \tag 2
$$
for $a<x<y<b$.  Notice that as $y\to a$ with $y>a$, the right-hand side
goes to $0$.  This proves that $T^+(a)=T$ exists.  
Likewise $T_-$ exists at every point.

Letting $x\to y$ and then $y\to a$ in the definition \thetag1 of $T^+(a)$, 
we immediately read off that $T_-(y)\to T^+(a)$ as $y\to a$.
If, however, we let $y\to x$ and then $x\to a$, 
we get that $T^+(x)\to T^+(a)$.  The convergence 
to $T_-(b)$ from the left is proved analogously.
\enddemo

\numbereddemo{Remark} It follows from \thetag2
that if $a\le x<y\le b$, then the 
angle between $T_{xy}$ and $T_{ab}$ is at most $2\kappa(a,b)$, 
where $\kappa(a,b)$ is the total
curvature of $\gamma|(a,b)$.  Consequently, the length of $\gamma(a,b)$
is at most $1/\cos(2\kappa(a,b))$ times $|\gamma(b)-\gamma(a)|$ provided
$\kappa(a,b)<\pi/4$.  
(This incidentally can be used to give another proof of Theorem 10.1.) 
It follows
that $T^+(a)$ is the right derivative of $\gamma$ with respect to arclength.
Likewise $T_-(a)$ is the corresponding left derivative.
\enddemo

\section{The collection of minimal M\"obius strips spanning
            a given curve}

\proclaim{Theorem}
Let $\GG$ be the set of smooth simple closed curves in $\RR^N$ with
total curvature $<4\pi${\rm .}  Let $\Cal N$ be the set of curves $\Gamma$ in
$\GG$ that bound minimal M{\rm \"{\it o}}bius strips with nonzero nullity{\rm .}
Then\/{\rm :}
\medbreak
\item{\rm (i)} the set $\Cal N$ is closed and nowhere dense{\rm ,} and
\smallbreak \item{\rm (ii)} each curve $\Gamma \in \GG\setminus \Cal N$
bounds only finitely many minimal M{\rm \"{\it o}}bius strips{\rm ,} exactly half of which have 
even index{\rm .}

\endproclaim

\demo{Proof}
Let $\MM$ be the set of all minimal M\"obius strips
bounded by curves in $\GG$, and let $\Pi$ be the forgetful 
map from $\MM$ to $\GG$:
$$ 
\align \noalign{\vskip-3pt}
&\Pi: \MM \to \GG, \\ \noalign{\vskip-3pt}
&\Pi: M \mapsto \partial M.
\endalign
$$
By Theorem 3.2, the elements of $\MM$
are all smooth and embedded.  Thus by  [Wh2, \S1.5 and \S3.3], 
$\MM$ has the structure of
a separable Banach manifold, 
and $\Pi$ is a smooth Fredholm map of Fredholm index $0$.

We claim that the map $\Pi$ is proper.  
To see this, let $M_i$ be a sequence
of surfaces in $\MM$ with $\Pi(M_i)=\partial M_i$ converging smoothly
to $\Gamma \in \GG$.  By the Gauss-Bonnet theorem, the surfaces
in $\MM$ have uniformly bounded total curvature.  Hence a subsequence 
of the $M_i$ will converge smoothly to a limit surface $M$ except at 
isolated singular points of $M$. 
(See parts (1) and (2) of Theorem 3 in \cite{Wh3}.
The theorem is stated only for $3$-manifolds, but the proof of parts 
(1) and (2) is valid for minimal surfaces in $N$-manifolds.)
Since $M$ has finite topological type, the singularities can only
be branch points.  (See \cite{O2, Theorem 1} for interior singularities;
a slight modification of the proof there establishes the result for
boundary singularities.)
But by Theorem 3.2, $M$ has no branch points.  
Hence $\Pi$ is a proper map.

By the Sard-Smale theorem \cite{Sm} 
or by  [Wh2, \S3.3(6)],
the set of critical values of $\Pi$ is 
a countable union of closed nowhere dense sets.  
By the properness of $\Pi$,
in fact this set must be closed and nowhere dense.  
By  [Wh2, \S3.3], the set of 
critical values is precisely $\Cal N$.  This proves \therosteritem{i}.

We also claim that the space $\GG$ is path connected.
The properness of $\Pi$ then implies  [Wh2, \S3.3]
that $\Pi$ has a mapping degree $d$ that, for any noncritical $\Gamma$, is 
equal to the difference between the numbers of even and odd index minimal
M\"obius strips bounded by $\Gamma$.  If $\Gamma$ is a planar curve,
then it bounds no minimal M\"obius strips.  Thus $d$ must be $0$.
This proves \therosteritem{ii}, given the path connectivity of $\GG$.

Thus it remains only to show that $\GG$ is path connected.
That is, we must show that a curve $\Gamma$ in $\GG$ can be deformed
through curves in $\GG$
to a plane convex curve (for example).
By Milnor, there is a linear function $L$ on $\RR^N$ such that
the restriction of $L$ to $\Gamma$ has
exactly one local maximum and one local minimum.  We may also assume that
the maximum is nondegenerate.
Let $[a,b]$ be the image of $\Gamma$ under $L$.
For $t\in [a,b)$, let $\Gamma(t)$ be the curve obtained by replacing the 
portion of $\Gamma$ where $L\le t$ by a straight line segment.
Note that each $\Gamma(t)$ is embedded.
Also, this truncation can only decrease total curvature, 
so $\Gamma(t)$ has total
curvature less than $4\pi$.  When $t=a$, $\Gamma(t)$ is the original curve,
and when $t$ is close to $b$, $\Gamma(t)$ is very nearly a convex plane curve
(namely, half of an ellipse).  Of course $\Gamma(t)$ is only piecewise smooth.
We leave it to the reader to modify the deformation so that each curve is
smooth and so that the final curve is planar and convex.
\enddemo

\bigbreak\centerline{\bf Remarks}

\bigbreak
{11.2}.
The theorem also holds, with the same proof, for surfaces of any 
other topological type except the disk type.  However, if the conjecture
in Section 5 is correct, then there are no such surfaces.

\bigbreak{11.3}.
The theorem also holds, with the same proof, for disk-type 
surfaces, except that in conclusion \therosteritem{ii} of the theorem,
the number of even index disks will be one more than the number of
odd index disks.  Of course when $N=3$, there is only one disk by Nitsche's 
theorem.

\bigbreak{11.4}.
For orientable surfaces, these assertions (11.2 and 11.3)
follow immediately from work of Tomi and Tromba (\cite{TT3}).  
Moreover, for $N\ge 4$, the assumption on total curvature is not needed.
Presumably the theory in \cite{TT3} could be 
modified for nonorientable surfaces.  

\bigbreak{11.5}.
Theorem 11.1 could presumably be proved by 
minimax theorems as in \cite{Sh2} or \cite{PR}. 
However, \cite{Sh2} assumes orientability,
and details of \cite{PR} have not yet appeared.

\AuthorRefNames [DHKW]
\references

[AA] \name{W. K. Allard} and \name{F. J. Almgren, Jr.},
The structure of stationary one dimensional varifolds
       with positive density,
{\it Invent.\ Math.\/}
{\bf 34} (1976), 83--97.

[All1] \name{W. K. Allard},
On the first variation of a varifold,
{\it Ann.\ of Math.\/}
{\bf 95} (1972), 417--491.

[All2]
\bibline,
On the first variation of a varifold: boundary behavior,
{\it Ann.\ of Math.\/}
{\bf 101} (1975), 418--446.

[Alt] \name{H. W. Alt},
Verzweigungspunkte von $H$-Fl\"achen, I,
{\it Math.\ Z.}
{\bf 127} (1972), 333--362; 
II,
{\it Math.\ Ann.\/}
{\bf 201} (1973), 33--55.

[AT] \name{F. J. Almgren, Jr.} and \name{W. P. Thurston},
Examples of unknotted curves which bound only surfaces of high
  genus within their convex hull, 
{\it Ann.\ of Math.\/}
{\bf 105} (1977), 527--538.

[AS] \name{F. J. Almgren, Jr.} and \name{L. Simon},
Existence of embedded solutions of Plateau's problem,
{\it Ann.\ Scuola Norm.\ Sup.\ Pisa}
{\bf 6} (1979), 447--495.

[Be] \name{F. Bernatzki},
The Plateau-Douglas problem for nonorientable minimal surfaces,
{\it Manu\-scripta Math.\/}
{\bf 79} (1993), 73--80.

[Boh]  \name{R. B\"ohme},
New results on the classical problem of Plateau.
       On the existence of many solutions,
 {\it S{\rm \'{\it e}}minaire Bourbaki},
{\it Ast\'erisque}
{\bf 92} (1981/82), 1--20.

 [Bor] \name{K. Borsuk},
Sur la courbure totale des courbes ferm\'ees,
{\it Annales Soc.\ Polonaise math.\/}
{\bf 20} (1947), 251--265.

[DHKW] \name{U. Dierkes, S. Hildebrandt, A. K\"uster}, and \name{O. Wohlrab},
{\it Minimal Surfaces},
Vols.\ I \& II,
 Springer-Verlag, New York, 1992.

[DW] \name{J. Drachman} and \name{B. White},
Soap films bounded by non-closed curves,
{\it J. Geom.\ Anal.\/}
{\bf 8} (1998), 239--250.

[Fa] \name{M. I. F\'ary},
Sur la courbure totale d'une courbe gauche faisant un n\oe ud,
{\it Bull.\ Soc.\ Math.\ France}
{\bf 77} (1949), 128--138.

[Fed] \name{H. Federer},
Some theorems on integral currents,
{\it Trans.\ Amer.\ Math.\ Soc.\/}
{\bf 117} (1965), 43--67.

[Fen] \name{W. Fenchel},
\"Uber Kr\"ummung und Windung geschlossener Raumkurven,
{\it Math.\ Ann.\/}
{\bf 101} (1929), 238--252.

[Gr] \name{M. Gromov},
Filling Riemannian manifolds,
{\it J. Differential Geom.\/}
{\bf 18} (1983), 1--147.

[Gu1] \name{R. Gulliver},
Regularity of minimizing surfaces of prescribed mean
       curvature,
{\it Ann.\ of Math.\/}
{\bf 97} (1973), 275--305.

[Gu2]
\bibline,
A minimal surface with an atypical boundary branch point,
in {\it Differential Geometry}, {\it Pitman Surveys Pure Appl.\ Math.\/} {\bf 52} (H. B. Lawson, Jr. and K. Tenenblat,
eds.)
Proc.\ of conference for M. do-Carmo, Longman/Wiley, Harlow/New York, 1991, 211--228.

[GuL] \name{R. Gulliver} and \name{F. Lesley},
On boundary branch points of minimizing surfaces,
{\it Arch.\ Rat.\ Mech.\ Anal.\/}
{\bf 52} (1973), 20--25.

[GuS] \name{R. Gulliver} and \name{J. Spruck},
On embedded minimal surfaces,
{\it Ann.\ of Math.\/}
{\bf 103} (1976), 331--347;
 correction 
{\it Ann.\ of Math.\/}
{\bf 109} (1979), 407--412.

[H] \name{S. Hildebrandt},
Boundary behavior of minimal surfaces,
{\it Arch.\ Rat.\ Mech.\ Anal.\/}
{\bf 35} (1969), 47--82.

[HS] \name{R. Hardt} and \name{L. Simon},
Boundary regularity and embedded solutions for the
       oriented Plateau problem,
{\it Ann.\ of Math.\/}
{\bf 110} (1979), 439--486.

[HT] \name{E. Heinz} and \name{F. Tomi},
Zu einem Satz von Hildebrandt
     \"uber das Randverhalten von Minimalfl\"achen,
{\it Math.\ Z.}
{\bf 111} (1969), 372--386.

[Hub] \name{J. H. Hubbard},
On the convex hull genus of space curves,
{\it Topology}
{\bf 19} (1980), 203--208.

[J] \name{J. Jost},
Conformal mappings and the Plateau-Douglas problem in Riemannian
        manifolds,
{\it J. Reine Angew.\ Math.\/}
{\bf 359} (1985), 37--54.

[K] \name{D. Kinderlehrer},
The boundary regularity of minimal surfaces,
{\it Ann.\ Scuola Norm.\ Sup.\ Pisa}
{\bf 23} (1969), 711--744.

[Ku] \name{N. Kuiper},
Geometry in curvature theory, 
in {\it Tight and Taut Submanifolds}\break
(T. E. Cecil and S.-s. Chern, eds.),
 {\it Math.\ Sci.\ Res.\ Inst.\ Publication} {\bf 32},
  Cambridge Univ.\ Press, Cambridge, 1997, 1--50.

[L] \name{F. D. Lesley},
Differentiability of minimal surfaces at the boundary,
{\it Pacific J. Math.\/}
{\bf 37} (1971), 123--139.

[Me] \name{W. H. Meeks, III},
Uniqueness theorems for minimal surfaces,
{\it Illinois J. Math.}\/
{\bf 25} (1981), 318--336.

[Mi] \name{J. Milnor},
On the total curvature of knots,
{\it Ann.\ of Math.\/}
{\bf 52} (1950), 248--257.

[Mo1] \name{F. Morgan},
A regularity theorem for minimizing hypersurfaces modulo $\nu$,
{\it Trans.\ Amer.\ Math.\ Soc.\/}
{\bf 297} (1986), 243--253.

[Mo2]
\bibline,
{\it Geometric Measure Theory\/{\rm :} A Beginner's Guide},
 Academic Press,  San Diego, CA, 1995.

[MW] \name{W. H. Meeks, III} and \name{B. White},
The space of minimal annuli bounded by an extremal pair
       of planar curves,
{\it Comm.\ Anal.\ Geom.\/}
{\bf 1} (1993), 415--437.

[MY1] \name{W. H. Meeks, III} and \name{S. T. Yau},
The classical Plateau problem and the topology of three-dimensional  
       manifolds.  The embedding of the solution given by Douglas-Morrey
       and an analytic proof of Dehn's lemma,
{\it Topology}
{\bf 21} (1982), 409--442.

[MY2]
\name{W. H. Meeks, III} and \name{S. T. Yau},
The existence of embedded minimal surfaces and the problem of 
          uniqueness,
{\it Math.\ Z.}
{\bf 179} (1982), 151--168.

[N1] \name{J. C. C. Nitsche},
MR 25 \#492, (review of Sasaki's paper [Sa]),
{\it Math.\ Rev.\/}
{\bf 25} (1963),
 104.

[N2]
\bibline,
The boundary behavior of minimal surfaces.  Kellog's theorem
       and Branch points on the boundary,
{\it Invent.\ Math.\/}
{\bf 8} (1969), 313--333; 
Concerning my paper on the boundary behavior of minimal surfaces,
{\it Invent.\ Math.\/}
{\bf 9} (1969/70), 270.

[N3]
\bibline,
A new uniqueness theorem for minimal surfaces,
{\it Arch.\ Rat.\ Mech.\ Anal.\/}
{\bf 52} (1973), 319--329.

[N4]
\bibline,
{\it Lectures on Minimal Surfaces}, Vol.\ I,
 Cambridge Univ.\ Press, Cambridge, 1989.

[O1] \name{R. Osserman},
A proof of the regularity everywhere of the classical
       solution to Plateau's problem,
{\it Ann.\ of Math.\/}
{\bf 91} (1970), 550--569.

[O2]
\bibline,
On Bers' theorem on isolated singularities,
{\it Indiana Univ.\ Math.\ J.}
{\bf 23} (1973), 337--342.

 [OS] \name{R. Osserman} and \name{M. Schiffer},
Doubly-connected minimal surfaces,
{\it Arch.\ Rational Mech.\ Anal.\/}
{\bf 58} (1975), 285--307.

[PR] \name{J. Pitts} and \name{H. Rubinstein},
Existence of minimal surfaces of bounded topological type
       in three-manifolds,
in  {\it Proc.\ Centre Math.\ Anal.\ Austral.\ Nat.\ Univ.\/}, Miniconference on geometry and partial differential equations
          (Canberra, 1985),
Vol.\ 10, Canberra, 1986, 163--176.

[Ra1] \name{T. Rad\'o},
On Plateau's problem
{\it Ann.\ of Math.\/}
{\bf 31} (1930), 457--469.

[Ra2]
\bibline,
{\it On the Problem of Plateau},
 Springer-Verlag, New York, 1971.

 [Sa] \name{S. Sasaki},
On the total curvature of a closed curve,
{\it Japanese J. Math.}\/
{\bf 29} (1959), 118--125.

[Sch] \name{R. Schneider},
A note on branch points of minimal surfaces,
{\it Proc.\ Amer.\ Math.\ Soc.\/}
{\bf 17} (1966), 1254--1257.

[Sh1] \name{M. Shiffman},
The Plateau problem for minimal surfaces of arbitrary topological
           structure,
{\it Amer.\ J. Math.\/}
{\bf 61} (1939), 853--882.

[Sh2]
\bibline, 
Unstable minimal surfaces with several boundaries,
{\it Ann.\ of Math.\/}
{\bf 43} (1942), 197--222.

[Si] \name{L. Simon},
{\it Lectures on Geometric Measure Theory}, 
{\it Proc.\ of the Centre for Math.\ Anal.\/}
{\bf 3},
  Australian National University Centre for Mathematical Analysis,
 Canberra, Australia, 1983, vii+272pp.

[Sm] \name{S. Smale},
An infinite dimensional version of Sard's theorem,
{\it Amer.\ J. Math.\/}
{\bf 87} (1965), 861--866.

[TT1] \name{F. Tomi} and \name{A. J. Tromba},
Extreme curves bound embedded minimal surfaces of the type of the
       disc,
{\it Math.\ Z.}
{\bf 158} (1978), 137--145.

[TT2]
\bibline,
Existence theorems for minimal surfaces of nonzero genus
       spanning a contour,
{\it Mem.\ Amer.\ Math.\ Soc.\/}
{\bf 71} (1988),
 iv+83 pp.

[TT3]
\bibline,
The index theorem for minimal surfaces of higher genus,
{\it Mem.\ Amer.\ Math.\ Soc.\/}
{\bf 117} (1995), vi+78 pp.

[War] \name{S. E. Warschawski},
Boundary derivatives of minimal surfaces,
{\it Arch.\ Rational Mech.\ Anal.\/}
{\bf 38} (1970), 241--256.

[Wh1] \name{B. White},
A regularity theorem for minimizing hypersurfaces mod $p$,
in {\it Geometric Measure Theory and the Calculus of Variations}, {\it Proc.\ Sympos.\ Pure Math.}\/ {\bf 44} 
(1986), Amer.\ Math.\ Soc., 413--427.

[Wh2]
\bibline,
The space of $m$-dimensional surfaces that are stationary for a
   parametric elliptic functional,
{\it Indiana Univ.\ Math.\ J.}
{\bf 36} (1987), 567--602.

[Wh3]
\name{B. White},
Curvature estimates and compactness 
       theorems in $3$-manifolds for surfaces that are 
       stationary for parametric elliptic functionals,
{\it Invent.\ Math.\/}
{\bf 88} (1987), 243--256.

[Wh4]
\name{B. White},
Half of Enneper's surface minimizes area,
in {\it Geometric Analysis and the Calculus of Variations for
        Stefan Hildebrandt} (J. Jost, ed.),
 International Press,  Cambridge, MA, 1996, 361--367.

[Wh5]
\bibline,
Classical area minimizing surfaces with real-analytic boundaries,
{\it Acta Math.\/}
{\bf 179} (1997), 295--305.

[Wi] \name{D. Wienholtz},
A method to exclude branch points of minimal surfaces,
{\it Calc.\ Var.\ Partial Differential Equations}
{\bf 7} (1998), 219--247.

\endreferences

\enddocument